\documentclass[12pt]{article}

\usepackage{amsfonts}

\newtheorem{thm}{Theorem}[section]
\newtheorem{lemma}[thm]{Lemma}
\newtheorem{cor}[thm]{Corollary}

\newtheorem{prop}[thm]{Proposition}
\newtheorem{conjecture}{Conjecture}

\newcommand{\proof
}{\par\medskip\noindent {\bf Proof.\ \ }}

\newcommand{\be}{\begin{equation}}
\newcommand{\ee}{\end{equation}}
\newcommand{\openbox}{\leavevmode
  \hbox to8pt{\hfil\vrule\vbox to6pt{\hrule width6pt\vfil\hrule}\vrule}}

\newcommand{\qed}{\hbox to5pt{ } \hfill \openbox\bigskip\medskip}

\newcommand{\Fp}{\mathbb F _p}

\newcommand{\cT}{\mbox{$\cal T$}}

\newcommand{\cF}{\mbox{$\cal F$}}
\newcommand{\cG}{\mbox{$\cal G$}}
\newcommand{\cH}{\mbox{$\cal H$}}

\newcommand{\cD}{\mbox{$\cal D$}}

\newcommand{\Z}{\mathbb Z}
\newcommand{\Q}{\mathbb Q}

\newcommand{\C}{\mathbb C}
\newcommand{\F}{\mathbb F}

\title{Balancing Sets of Vectors}
\author{G\'abor Heged\H{u}s
\\{\normalsize  \'Obuda  University}
}

\begin{document}

\maketitle

\begin{abstract}
Let $n$ be an arbitrary integer, let $p$ be a prime factor of $n$. Denote
by $\omega_1$ the $p^{th}$ primitive unity root, $\omega_1:=e^{\frac{2\pi i}{p}}$.

Define $\omega_i:=\omega_1^i$ for $0\leq i\leq p-1$ and $B:=\{1,\omega_1,\ldots,\omega_{p-1}\}^n\subseteq \C^n$.

Denote by $K(n,p)$ the minimum $k$ for which there exist vectors $v_1,\ldots,v_k\in B$ such that for any vector $w\in B$, there is an $i$, $1\leq i\leq k$, such that $v_i\cdot w=0$, where $v\cdot w$ is the usual scalar product of $v$ and $w$.

Gr\"obner basis methods and linear algebra proof gives the lower bound $K(n,p)\geq n(p-1)$.

Galvin posed the following problem: Let $m=m(n)$ denote the minimal integer such that there exists subsets $A_1,\ldots,A_m$ of $\{1,\ldots,4n\}$ with $|A_i|=2n$ for each $1\leq i\leq n$, such that for any subset $B\subseteq [4n]$ with $2n$ elements there is at least one $i$, $1\leq i\leq m$, with $A_i\cap B$ having $n$ elements. We obtain here the result $m(p)\geq p$ in the case of $p>3$ primes.
\end{abstract}

\section{Introduction}

First we introduce some notations.

Let $n$ be an arbitrary positive integer and consider a $p$ prime factor of $n$. Denote
by $\omega_1$ the $p^{th}$ primitive unity root, i.e., let $\omega_1:=e^{\frac{2\pi i}{p}}$. Define $\omega_i:=\omega_1^i$ for each $1\leq i\leq p-1$.

Throughout the paper $\F$ denotes a fixed field. As usual,
$\F[x_1, \ldots, x_n]$ denotes  the
ring of polynomials in variables $x_1, \ldots, x_n$ over $\F$.
We use also the shorter notation $S=\F[x_1,\ldots,x_n]$.

Let $[n]$ stand for the set $\{1,2,
\ldots, n\}$. For an integer $0\leq d\leq n$ we denote by
${[n] \choose d}$ the family of all  $d$ element subsets of $[n]$.

Let $v_F\in \{0,1\}^n$ denote the
characteristic vector of a set
$F \subseteq [n]$.
For a family of subsets $\cF \subseteq 2^{[n]}$, define
$$
V(\cF):= \{v_F : F \in \cF\} \subseteq \{0,1\}^n \subseteq \F^n
$$
the set of characteristic vectors of the family $\cF$.

Let $R(n,d)$ denote the minimal $k$ for which there exist vectors $v_1,\ldots,v_k\in \{-1,1\}^n$ such that for any vector $w\in \{-1,1\}^n$ there is an $i$, $1\leq i\leq k$ such that $|v_i\cdot w|\leq d$, where $v\cdot w$ denotes the usual inner product of two vectors. Since $v\cdot w\equiv n \pmod 2$ for any two vectors $v,w\in \{-1,1\}^n$, $R(n,0)$ is defined only for even $n$, while $R(n,d)$ for $d\geq 1$ is well--defined for all $n$. A simple construction of Knuth \cite{K} shows that $R(n,d)\leq \lceil n/(d+1)\rceil$ for $n\equiv d \pmod 2$, where $\lceil x \rceil$ denotes the least integer which is at least $x$. In \cite{ABCO} Alon, Bergmann, Coppersmith and Odlyzko showed that this construction is optimal. They used in their proof only elementary linear algebra.

It is possible to generalize this problem and consider balancing families of vectors whose components are $p^{th}$ root of unity for some fixed $p$. Our main result is the following:

\begin{thm} \label{main2}
Let $n$ be an arbitrary integer, let $p$ be a prime factor of $n$. Denote
by $\omega_1$ the $p^{th}$ primitive unity root, $\omega_1:=e^{\frac{2\pi i}{p}}$.

Define $\omega_i:=\omega_1^i$ for $0\leq i\leq p-1$ and $B:=\{1,\omega_1,\ldots,\omega_{p-1}\}^n\subseteq \C^n$.

Denote by $K(n,p)$ the minimum $k$ for which there exist vectors $v_1,\ldots,v_k\in B$ such that for any vector $w\in B$, there is an $i$, $1\leq i\leq k$, such that $v_i\cdot w=0$, i.e., $v$ is orthogonal with respect to the usual scalar product to $w$. Then $K(n,p)\geq n(p-1)$. \qed
\end{thm}

We suggest the following plausible conjecture:
\begin{conjecture}
Let $n$ be an arbitrary integer, let $p$ be a prime factor of $n$. Then $K(n,p)=n(p-1)$.
\end{conjecture}

We can rephrase the previous balancing vector problem in term of an extremal combinatorial problem for subsets of a set, with an $n$-dimensional vector $u=(u_1,\ldots,u_n)\in \{-1,1\}^n$ corresponding a subset $A$ of $\{1,2\ldots,n\}$ with $j\in A$ iff $u_j=1$. Galvin posed a similar problem in this setting. He asked for a determination of the minimal integer $m=m(n)$ such that there exist subsets $A_1,\ldots,A_m$ of $\{1,\ldots,4n\}$, $|A_i|=2n$ for each $1\leq i\leq m$, with the following property: for any subset $B\subseteq [4n]$ with $2n$ elements there is at least one $i$, $1\leq i\leq m$, with $A_i\cap B$ having $n$ elements.

Galvin noticed that if one defines $A_i=\{i,i+1,\ldots , i+2n-1\}$ for $1\leq i\leq 2n$, then it is easy to verify that these $A_i$ have the right property, so $m(n)\leq 2n$.

We obtain the following Theorem with an other application of Gr\"obner basis methods and linear algebra.
\begin{thm} \label{main}
Let $p>3$ be a prime. Then $m(p)\geq p$.
\end{thm}

The organisation of this article is the following:

In Section 2 we define Gr\"obner bases and standard monomials in polynomial rings. In Section 3 we prove our main method giving a general lower bound for the degree of a polynomial via standard monomials. In Section 4 we determine the standard monomials of combinatorially interesting finite subsets. In Section 5 we prove our main results.

\medskip

\section{Gr\"obner bases and standard monomials}

We recall now some basic facts
concerning Gr\"obner bases in polynomial rings.
A total order $\prec$ on the monomials $\mbox{Mon}$ of the polynomial ring is a {\em term order}, if 1 is the
minimal element of $\prec$, and $uw\prec vw$ holds for any monomials
$u,v,w$ with $u\prec v$. We define
now two interesting term orders: the
lexicographic (lex) and the
deglex term orders. Let $u=x_1^{i_1}x_2^{i_2}\cdots x_n^{i_n}$ and
$v=x_1^{j_1}x_2^{j_2}\cdots x_n^{j_n}$ be two monomials. Then $u$ is smaller
than $v$ with respect to lex ($u\prec_{lex} v$ in notation) iff $i_k<j_k$ holds for the smallest index $k$ such
that $i_k\not=j_k$. Similarly, $u$ is smaller
than $v$ with respect to deglex ($u\prec_{deg} v$ in notation) iff either
$\deg u< \deg v$, or $\deg u=\deg v$ and $u\prec_{lex} v$. Note that we have $x_n\prec x_{n-1}\prec\ldots \prec
x_1$, for both lex and deglex.
Clearly the deglex $\prec_{deg}$ order is a degree-compatible
term order (this means that $\deg u <\deg v $ implies $u\prec v$).

The {\em leading monomial} ${\rm lm}(f)$
of a nonzero polynomial $f\in S$ is the largest
(with respect to $\prec$) monomial
which appears with nonzero coefficient in $f$ when written as a linear
combination of different monomials.
The {\em initial ideal} $\mbox{in}(I)$ of an ideal $I$ is
the set of all leading monomials $\mbox{lm}(f)$:
$\mbox{in}(I)=\{\mbox{lm}(f):~ f\in I\}$.

Let $I$ be an ideal of $S$. A finite subset $\cG\subseteq I$ is a {\it
Gr\"obner basis} of $I$ if for every $f\in I$ there exists a $g\in \cG$ such
that ${\rm lm}(g)$ divides ${\rm lm}(f)$. In other words, the leading
monomials of the polynomials from $\cG$ generate the semigroup ideal of monomials
$\mbox{in}(I)$. It follows from the fact that $\prec$ is a well founded
order, that $\cG$ is
actually a basis of $I$, i.e.,  $\cG$ generates $I$ as an ideal of $S$. We can prove easily (cf. \cite[Chapter 1, Corollary
3.12]{CCS} or \cite[Corollary 1.6.5, Theorem 1.9.1]{AL}) that every
nonzero ideal $I$ of $S$ has a Gr\"obner basis.

A Gr\"obner basis $\{ g_1,\ldots ,g_m\}$ of $I$ is {\em reduced}
if the coefficient of ${\rm lm}(g_i)$ is 1, and no nonzero monomial
in $g_i$ is divisible by any ${\rm lm}(g_j)$, $j\not= i$.
By a theorem of Buchberger (\cite[Theorem 1.8.7]{AL}) a nonzero ideal
has a unique reduced Gr\"obner basis.

A monomial $w\in S$ is called a {\it standard monomial for $I$} if
it is not a leading monomial of any $f\in I$. Let ${\rm Sm}(\prec,I)$ stand
for the
set of all standard monomials of $I$ with respect to the term-order
$\prec$ over $\F$.
Using the definition and existence of
Gr\"obner bases (see \cite[Chapter 1, Section 4]{CCS}) we can prove easily that for a
nonzero ideal $I$  the set
${\rm Sm}(\prec,I)$ is a basis of the $\F$-vector-space $S/I$. More
precisely,
every $g\in S$ can be written uniquely as $g=h+f$ where $f\in I$ and
$h$ is a unique $\F$-linear combination of monomials from ${\rm Sm}(\prec,I)$. We say that the polynomial $h$
is the {\em reduction} of the polynomial $g$ via the Gr\"obner basis $\cG$ of the ideal $I$.

In general {\em reduction} means that we repeatedly replace monomials
in $f$ by smaller ones (with respect to $\prec$). The procedure is the following: if monomial $w$ occurs in $f$ and ${\rm lm}(g)$ divides $w$ for some
$g\in {\cal G}$, then
we replace $w$ in $f$ with $w-\frac{w}{{\rm lm}(g)}\cdot g$. Obviously the
monomials in $w-\frac{w}{{\rm lm}(g)}\cdot g$ are $\prec$-smaller than $w$.

For $\cF\subseteq \F^n$, $\cF\not=\emptyset $ we put
$$
\mbox{Sm}(\prec,\cF):={\rm Sm}(\prec,I(\cF)),
$$
where
$$
I(\cF):=\{f\in S:~f(v)=0 \mbox{ whenever } v\in \cF\}.
$$

It is clear that $\mbox{Sm}(\prec,\cF)$ is downward closed, i.e., if $u\in\mbox{Sm}(\prec,\cF)$ and $w$ divides $u$, then $w\in\mbox{Sm}(\prec,\cF)$.

Also, the standard monomials for $I(\cF)$ form a basis of the functions
from $\cF$ to $\F$, hence
\begin{equation} \label{egyenlo}
|\mbox{Sm}(\prec,\cF)|=|\cF|.
\end{equation}

\section{The method}

First we prove a general condition which gives a lower bound for the degree of a polynomial.

\begin{thm} \label{method}
Let $\F$ be an arbitrary field and $P(x_1,\ldots ,x_n) \in \F[x_1,\ldots,x_n]$ be an arbitrary polynomial.

Let $\cF \subseteq \F^n$ denote an arbitrary finite subset of the affine space such that $\cF\neq \F^n$ and let $\underline{h}\in \F^n\setminus \cF$. We put $\cT:= \cF\cup \{\underline{h}\}$.

Suppose that $P(\underline{h})\neq 0$ and $P(\underline{f})=0$ for each $\underline{f}\in \cF$. Let
$$
y\in \mbox{Sm}(\prec_{deg},\cT)\setminus \mbox{Sm}(\prec_{deg},\cF).
$$
Then $deg(P)\geq deg(y)$.
\end{thm}

\proof

Write $\cG$ for the deglex Gr\"obner basis of the ideal $I(\cT)$. We denote by $\overline{P}$ the reduction of $P$ via the Gr\"obner basis $\cG$. Then $deg(\overline{P})\leq deg(P)$, because in the process of reduction we replaced each monomial of $P$ with such monomials which have smaller degree. Clearly $\overline{P}(\underline{h})=P(\underline{h})\neq 0$, $\overline{P}(\underline{f})=P(\underline{f})=0$ for each $\underline{f}\in \cF$, because we reduced $P$ with such polynomials which vanish on $\cT$.

We can expand $\overline{P}$ into the unique form
\begin{equation} \label{expand}
\overline{P}=\sum_{m\in \mbox{Sm}(\prec_{deg},\cT)} \alpha_m\cdot m,
\end{equation}
where $\alpha_m\in \F$. It is enough to prove that $\alpha_y\neq 0$, namely then $deg(\overline{P})\geq deg(y)$.

Suppose indirectly, that $\alpha_y=0$. Since $\cF\subseteq \cT$, thus $\mbox{Sm}(\prec_{deg},\cF)\subseteq \mbox{Sm}(\prec_{deg},\cT)$ and $\mbox{Sm}(\prec_{deg},\cT)\setminus \mbox{Sm}(\prec_{deg},\cF)=\{y\}$. Therefore the equation (\ref{expand}) yields to the following expansion:
\begin{equation}
\overline{P}=\sum_{m\in \mbox{Sm}(\prec_{deg},\cF)} \alpha_m\cdot m,
\end{equation}
and since $\overline{P}(\underline{f})=0$ for each $\underline{f}\in \cF$, hence $\alpha_m=0$ for each $m\in \mbox{Sm}(\prec_{deg},\cF)$. But then $\overline{P}\equiv 0$ as functions mapping $\cT$ to $\F$, which gives a contradiction with $\overline{P}(\underline{h})\neq 0$. \qed

J. Farr and S. Gao proved in Lemma 2.2 of \cite{FG} the following.
\begin{lemma} \label{alap}
Suppose that $\cG=\{g_1,\ldots,g_s\}$ is a reduced Gr\"obner basis for the ideal $I(\cF)$, where $\cF\subseteq \F^n$ is a finite set of points. For a point $\underline{h}=(a_1,\ldots,a_n)\notin \cF$, let $g_i$ denote the polynomial in  $\cG$ with smallest leading term such that $g_i(\underline{h})\neq 0$, and define
\begin{equation}
\overline{g_j}:=g_j-\frac{g_j(\underline{h})}{g_i(\underline{h})}\cdot g_i,\mbox{   }j\neq i,\mbox{  and}
\end{equation}
\begin{equation}
g_{ik}:=(x_k-a_k)\cdot g_i,\mbox{    }1\leq k\leq n.
\end{equation}
Then
\begin{equation}
\overline{\cG}=\{\overline{g_1},\ldots,\overline{g_{i-1}},\overline{g_{i+1}},\ldots,\overline{g_s},g_{i1},\ldots,g_{in}\}
\end{equation}
constitutes a Gr\"obner basis for the ideal $I(\cF\cup \{\underline{h}\})$.
\end{lemma}

\begin{cor} \label{basic}
Let $\F$ be an arbitrary field and $\prec$ be an arbitrary term order on the monomials of $\F[x_1,\ldots,x_n]$. Let $\cF\subseteq \F^n$ stand for an arbitrary finite subset. Let $\underline{h}\in \F^n \setminus \cF$ be an arbitrary vector and define $\cT:= \cF\cup \{\underline{h}\}$.

Let $\cG=\{g_1,\ldots,g_s\}\subseteq \F[x_1,\ldots,x_n]$ stand for the reduced Gr\"obner basis of the ideal $I(\cF)$ with respect to the term order $\prec$.

Suppose that $m_1\prec \ldots \prec m_k$, where $m_i:=\mbox{lm}_{\prec}(g_i)$. Consider
$$
i:=\min \{j\in [k]:~ g_j(\underline{h})\neq 0\}.
$$
Then $\mbox{Sm}(\cT,\prec)=\mbox{Sm}(\cF,\prec)\cup \{m_i\}$.
\end{cor}

\proof

This Corollary is obvious from Lemma \ref{alap}. Namely
$$
|\mbox{Sm}(\prec,\cT)|=|\mbox{Sm}(\prec,\cF))|+1,
$$
therefore it is enough to prove that $m_i\in \mbox{Sm}(\prec,\cT)$.

Indirectly, suppose that $m_i\notin \mbox{Sm}(\prec,\cT)$. This means that there exists a polynomial $g\in \overline{\cG}$ such that $\mbox{lm}(g)$ divides $m_i$. Clearly if $j<i$, then $\mbox{lm}(\overline{g_j})=\mbox{lm}(g_j)=m_j$. Similarly, if $j>i$, then $\mbox{lm}(\overline{g_j})=\max(\mbox{lm}(g_j),\mbox{lm}(g_i))=\mbox{lm}(g_j)=m_j$.

Since $\cG$ was a reduced Gr\"obner basis of the ideal $I(\cF)$, hence $\mbox{lm}(g_j)=m_j$ does not divide $m_i$ for each $j\neq i$. Since
$$
\mbox{lm}(g_{il})=x_l\cdot \mbox{lm}(g_i)=x_l\cdot m_i,
$$
thus $\mbox{lm}(g_{il})$ does not divide also $m_i$ for each $1\leq k\leq n$, which gives a contradiction.
\qed

\begin{cor}
Let $\F$ be an arbitrary field and $\prec$ be an arbitrary term order on the monomials of $\F[x_1,\ldots,x_n]$. Let $\cF\subseteq \F^n$ stand for an arbitrary finite subset. Let $\underline{h}\in \F^n \setminus \cF$ be an arbitrary vector and put $\cT:= \cF\cup \{\underline{h}\}$.

Let $\cG=\{g_1,\ldots,g_s\}\subseteq \F[x_1,\ldots,x_n]$ stand for the reduced Gr\"obner basis of the ideal $I(\cF)$ with respect to the term order $\prec$.

Let $\chi_{\underline{h}}:\cT \to \F$ denote the characteristic function of $\underline{h}$, i.e., $\chi_{\underline{h}}(\underline{h})=1$ and $\chi_{\underline{h}}(\underline{f})=0$ for each $\underline{f}\in \cF$. Then
\begin{equation}
\chi_{\underline{h}}\equiv \frac{1}{g_i(\underline{h})}\cdot g_i
\end{equation}
gives an expansion of $\chi_{\underline{h}}$ into the unique linear combination of standard monomials of the ideal $I(\cT)$.
\end{cor} \qed

\section{Standard monomials}

Let $n$ be an arbitrary integer, let $p$ be a prime factor of $n$. Denote
by $\omega_1$ the $p^{th}$ primitive unity root, i.e., let $\omega_1:=e^{\frac{2\pi i}{p}}$. Define $\omega_i:=\omega_1^i$ for each $1\leq i\leq p-1$.
Write $B:=\{1,\omega_1,\ldots,\omega_{p-1}\}^n\subseteq \C^n$ and
$$
D:=\{x^u=x_1^{u_1}\cdot \ldots \cdot x_n^{u_n}:~ 0\leq u_i\leq p-1 \mbox{ for each }1\leq i\leq n\}.
$$

Let $B_j:=\{(t_1,\ldots,t_n)\in B:~ t_1\cdot \ldots \cdot t_n=\omega_j\}$ for each $0\leq j\leq p-1$. First we characterize the standard monomials and the reduced Gr\"obner basis of the ideal $I(B_0)\subseteq \C[x_1,\ldots,x_n]$ with respect to any $\prec$ term order.

Consider the following equivalence relation $\equiv$ on $D$:

let the monomials $x^u=x_1^{u_1}\cdot\ldots \cdot x_n^{u_n}$ and $x^v=x_1^{v_1}\cdot \ldots \cdot x_n^{v_n}$ be equivalent via $\equiv$ iff there exists a $k$, $0\leq k\leq p-1$ such that
$u_i+k\equiv v_i \pmod p$ for each $1\leq i\leq n$.

Denote by $D/\equiv$ the set of equivalence classes of $D$ with respect to $\equiv$ and write $[a]:= \{b\in D:~ b\equiv a\}$ for the equivalence class of $a\in D$. It is easy to verify that $|[a]|=p$ for each equivalence classes $[a]\in D/\equiv$, therefore $|D/\equiv|=p^{n-1}$.

Let $\prec$ be a fixed term order on the monomials of $\C[x_1,\ldots ,x_n]$.  Denote by $\mbox{Min}(\prec)$ the set of monomials $u$ of $D$ such that there exists an equivalence class $[a]\in D/\equiv$ for which $u$ is the minimal element of $[a]$ with respect to the term order $\prec$. Clearly $|\mbox{Min}(\prec)|=p^{n-1}$.

\begin{lemma} \label{vanish}
Let $[b]\in D/\equiv$ be an arbitrary equivalence class. Let $a$ denote the minimal element of $[b]$ with respect to the term order $\prec$ and suppose that $b\neq a$. Then there exists an $0\leq m\leq p-1$ such that the polynomial $b-\omega_m\cdot a\in I(B_j)$, where $0\leq j\leq p-1$.
\end{lemma}

\proof

Let $b:=x^u$ and $a:=x^v$. By the definition of the equivalence relation $\equiv$, $x^u\equiv x^v$ iff there exists a $k$, $0\leq k\leq p-1$ such that $u_i+k\equiv v_i \pmod p$ for each $1\leq i\leq n$. This means that $x^u$ is the reduction of the monomial $x^v\cdot (x_1\cdot\ldots \cdot x_n)^k$ via the polynomials $x_l^p-1$, where $1\leq l\leq n$. Since $B_j\subseteq B$ and $x_l^p-1\in I(B)$ for each $1\leq l\leq n$, hence $x_l^p-1\in I(B_j)$, and $x_1\cdot\ldots \cdot x_n-\omega_j\in I(B_j)$ by the definition of $B_j$, therefore $x^u(b)=\omega^{k\cdot j}x^v(b)$ for each $b\in B_j$. This gives that $x^u-\omega^{k\cdot j}x^v\in I(B_j)$.
\qed

\begin{prop} \label{stdmon}
Let $\prec$ be an arbitrary term order on the monomials of $\C[x_1,\ldots,x_n]$. Then  $\mbox{Sm}(\prec,B_i)=Min(\prec)$ for each $0\leq i\leq p-1$.
\end{prop}

\proof
Clearly
$$
|\mbox{Sm}(\prec,B_i)|=|B_i|=p^{n-1}=|\mbox{Min}(\prec)|.
$$

If $b=x_1^{u_1}\cdot\ldots \cdot x_n^{u_n}\notin D$, then $b\in \mbox{in}(I(B_i))$. Namely there exists an index $t$, $1\leq t\leq n$ such that $u_t\geq p$. Let $c$ denote the reduction of $b$ via $x_t^p-1$. Clearly $c\neq b$, and $b-c\in I(B)\subseteq I(B_i)$.

Therefore it is enough to show that for each $b\in D\setminus \mbox{Min}(\prec)$ there exists a polynomial $g_b\in I(B_i)$ such that $\mbox{lm}_{\prec}(g_b)=b$. Consider the equivalence class $[b] \in D/\equiv$ and let $a\in D$ denote the minimal element of this equivalence class with respect to the term order $\prec$. By Lemma \ref{vanish} there exists an $0\leq m\leq p-1$ such that the polynomial $g_{b}:=b-\omega_m\cdot a\in I(B_i)$. Since $b\notin \mbox{Min}(\prec)$, therefore $b\neq a$. It follows from the definition of $a$  that $\mbox{lm}_{\prec}(g_b)=b$.  \qed

\begin{thm} \label{Grobner}
Let $\prec$ be an arbitrary term order on the monomials of $\C[x_1,\ldots,x_n]$. Then the following set of polynomials constitute a reduced Gr\"obner basis  of the ideal $I(B_0)$ with respect to the term order $\prec$:
$$
\cG:=\{b-a:~ a \mbox{ is the minimal element of }[b],\ b\neq a,\ [b]\in D/\equiv\}
$$
$$
\cup\{x_i^p-1:~ 1\leq i\leq n\}.
$$
\end{thm}

\proof
To show that ${\cal G}$ is a Gr\"obner basis of $I(B_0)$ it is enough to prove that $\cG\subseteq I(B_0)$ and there exists a polynomial $g\in \cG$ for each $f\in I(B_0)$ such that $\mbox{lm}(g)$ divides $\mbox{lm}(f)$.

The containment $\cG\subseteq I(B_0)$ follows from Lemma \ref{vanish}.

Let $f\in I(B_0)$ be an arbitrary polynomial. Then $b:=\mbox{lm}(f)\notin \mbox{Sm}(\prec,B_0)=\mbox{Min}(\prec)$ by Proposition \ref{stdmon}. If $b=x_1^{u_1}\cdot\ldots\cdot  x_n^{u_n}\notin D$, then there exists an index $1\leq i\leq n$ such that $u_i\geq p$. Then clearly $\mbox{lm}(x_i^p-1)=x_i^p$ divides $u$.

If $b\in D\setminus Min(\prec)$, then let $a$ denote the minimal element of the equivalence class $[b]$. Then $g_b:=b-a$ gives our statement.

It is obvious from Proposition \ref{stdmon} that the leading terms of the polynomials in $\cG$ constitute the minimal generating set of the initial ideal of $I(B_0)$. Reducedness follows from the fact that all non-leading monomials in these polynomials are actually standard monomials for $I(B_0)$ by Proposition \ref{stdmon}.

\qed

We prove the following easy consequence of the characterization of standard monomials:

\begin{prop} \label{degmon}
Let $\prec$ be an arbitrary degree-compatible term order. Then
\begin{equation}
\{x^u\in D:~ deg(x^u)<\frac{n(p-1)}{p}\}\subseteq \mbox{Sm}(B_t,\prec)
\end{equation}
for every $0\leq t\leq p-1$.
\end{prop}

\proof

Let $0\leq t\leq p-1$ be fixed and let $x^u=b_0\in D$ be an arbitrary monomial and we denote by $b_k$ the reduction of $x^u\cdot (x_1\cdot \ldots \cdot x_n)^k/\omega^{k\cdot t}$ via the equations $x_i^p-1$, $1\leq i\leq n$, for each $0\leq k\leq p-1$. Suppose that $\mbox{deg}(b_0)< \frac{n(p-1)}{p}$. Then by Proposition \ref{stdmon} it is enough to prove that
\begin{equation} \label{csill}
\mbox{deg}(b_i)>\mbox{deg}(b_0)
\end{equation}
for each $1\leq i\leq p-1$, because $\prec$ was a degree-compatible term order, thus (\ref{csill}) means that $b_0$ is the minimal element of the equivalence class $[b_0]$.

We may suppose without lost of generality that
$$
b_0=x_1^{p-1}\cdot \ldots \cdot x_{\lambda_1}^{p-1}x_{\lambda_1+1}^{p-2}\cdot\ldots \cdot x_{\lambda_1+\lambda_2}^{p-2}\cdot \ldots \cdot x_{\lambda_1+\ldots +\lambda_{p-2}+1}^1\cdot \ldots \cdot x_{\lambda_1+\ldots +\lambda_{p-2}+\lambda_{p-1}}^1,
$$
where $n=\lambda_1+\ldots +\lambda_p$.

Then
\begin{equation} \label{fok}
\mbox{deg}(b_0)=\sum_{j=1}^{p-1} (p-j)\cdot\lambda_j<\frac{n(p-1)}{p}.
\end{equation}

It is easy to verify from the definition of $b_i$ that
$$
b_i=x_1^{i-1}\cdots  x_{\lambda_1}^{i-1}x_{\lambda_1+1}^{i-2}\cdots x_{\lambda_1+\lambda_2}^{i-2}\cdots x_{\lambda_1+\ldots +\lambda_i+1}^{p-1}\cdots x_{\lambda_1+\ldots +\lambda_{i+1}}^{p-1}\cdots x_{\lambda_1+\ldots +\lambda_{p-1}+1}^{i}\cdots x_{\lambda_1+\ldots +\lambda_p}^{i}.
$$
Then
$$
\mbox{deg}(b_i)=\sum_{j=1}^{i-1} (i-j)\cdot\lambda_j+\sum_{j=1}^{p-i} (p-j)\cdot\lambda_{i+j}.
$$
Therefore it is enough to prove that
$$
\sum_{j=1}^{p-1} (p-j)\cdot\lambda_j< \sum_{j=1}^{i-1} (i-j)\cdot\lambda_j+\sum_{j=1}^{p-i} (p-j)\cdot\lambda_{i+j}.
$$

This inequality is equivalent with
\begin{equation} \label{fok2}
(p-i)(\sum_{j=1}^i\lambda_j)< i(\sum_{j=i+1}^p\lambda_j)
\end{equation}
for each $1\leq i\leq p-1$.

It is easy to verify that the inequality (\ref{fok2}) is equivalent with
\begin{equation} \label{kocka}
(\sum_{j=1}^i\lambda_j)(p(p-i)-(p-1))<(\sum_{j=i+1}^p\lambda_j)\frac{i}{p-i}(p(p-i)-(p-1)).
\end{equation}

But $n=\lambda_1+\ldots +\lambda_p$, hence from (\ref{fok}) we get
\begin{equation} \label{fok3}
\sum_{j=1}^{p-1} (p-j)\cdot\lambda_j< \frac{p-1}{p}(\sum_{j=1}^p \lambda_j).
\end{equation}

After some rearrangement of the inequality (\ref{fok3}) we find that
\begin{equation} \label{fok4}
\lambda_1(p-1)^2+\ldots +\lambda_i(p(p-i)-(p-1))< \lambda_{i+1}(ip-(p-1)^2)+\ldots +\lambda_{p-1}(-1)+(p-1)\lambda_p.
\end{equation}

Now it is easy to verify that
\begin{equation} \label{fok5}
(\sum_{j=1}^i\lambda_j)(p(p-i)-(p-1))\leq \lambda_1(p-1)^2+\ldots +\lambda_i(p(p-i)-(p-1)).
\end{equation}

From (\ref{fok4}) and (\ref{fok5}) we conclude that
\begin{equation} \label{kocka2}
(\sum_{j=1}^i\lambda_j)(p(p-i)-(p-1))<\lambda_{i+1}(ip-(p-1)^2)+\ldots +(-1)\lambda_{p-1}+\lambda_p (p-1).
\end{equation}
But since
$$
(p-i)(p-1)\leq i(p(p-i)-(p-1))
$$
for each $1\leq i\leq p-1$ and $jp-(p-1)^2<0$ for each $i+1\leq j\leq p-1$, hence we get
\begin{equation} \label{kocka3}
\lambda_{i+1}(ip-(p-1)^2)+\ldots +\lambda_{p-1}(-1)+\lambda_p(p-1) < \frac{i}{p-i}(p(p-i)-(p-1))(\sum_{j=i+1}^{p-1}\lambda_j+\lambda_p)
\end{equation}
and the inequality (\ref{kocka}) follows from (\ref{kocka2}) and (\ref{kocka3}).

\qed

\begin{cor} \label{becs}
Let $0\leq t\leq p-1$ be an integer and let $\underline{q}\in B\setminus B_t$ be an arbitrary vector. Define $Q:=B_t\cup \{\underline{q}\}$ and consider $y\in \mbox{Sm}(\prec_{deg},Q)\setminus \mbox{Sm}(\prec_{deg},B_t)$. Then $deg(y)\geq \frac{n(p-1)}{p}$.
\end{cor}

\proof

Clearly $\mbox{Sm}(Q,\prec_{deg})\subseteq D$, hence $y\in D\setminus \mbox{Sm}(B_t,\prec_{deg})$. Since by Proposition \ref{degmon} $\{x^u\in D:~ deg(x^u)<\frac{n(p-1)}{p}\}\subseteq \mbox{Sm}(B_0,\prec)$, this means that $D\setminus \mbox{Sm}(B_t,\prec_{deg})\subseteq D\setminus \{x^u\in D:~ deg(x^u)<\frac{n(p-1)}{p}\}=\{x^u\in D:~ deg(x^u)\geq \frac{n(p-1)}{p}\}$.
\qed

\medskip

Now we characterize the standard monomials and the reduced Gr\"obner basis of the ideal $I(V{[4p] \choose 2p})\subseteq \Fp[x_1,\ldots ,x_{4p}]$, where $p$ is an arbitrary prime.

Let $n,t$ be integers such that $0<t\leq n/2$. We define
$\cH _t$ as the set of those subsets $\{s_1<s_2<\cdots <s_t\}$
of $[n]$ for which $t$ is the smallest index $j$ with  $s_j<2j$.

We get $\cH _1=\{\{1\}\}$,
$\cH _2=\{\{2,3\}\}$, and $\cH _3=\{\{2,4,5\},\{3,4,5\}\}$. It is clear
that if $\{s_1< \ldots <s_t\}\in \cH _t$, then $s_t=2t-1$, moreover
$s_{t-1}=2t-2$ if $t>1$.

For a subset $J\subseteq [n]$ and an integer $0\leq i\leq |J|$ we
denote by $\sigma_{J,i}$ the $i^{th}$ elementary symmetric polynomial of
the variables $x_j$, $j\in J$:
$$ \sigma _{J,i}:=\sum _{T\subseteq J, |T|=i} x_T~ \in \Z[x_1,\ldots
,x_n]. $$
Specifically, we have $\sigma _{J,0}=1$.

Now let $0<t\leq n/2$, $0\leq d\leq n$ and $H\in \cH_t$. Put $H'=H\cup
\{2t,2t+1,\ldots ,n\}\subseteq [n]$. We write
$$ f_{H,d}=f_{H,d}(x_1,\ldots, x_n):=\sum _{k=0}^t (-1)^{t-k}{d-k
\choose t-k}\sigma_{H',k}. $$
For example, we have $f_{\{1\},d}=x_1+x_2+\cdots +x_n-d$, and
$$ f_{\{2,3\},d}=\sigma_{U,2}-(d-1)\sigma_{U,1}+{d \choose 2}, $$
where $U=\{2,3,\ldots ,n\}$.

Let $\cD_d$ denote the collection of subsets
$x_U$, where $U=\{u_1<\ldots <u_{d+1}\}$ and $u_j\geq 2j$ holds for
$j=1,\ldots ,d$.

The following statement was proved in \cite{HR}.

\begin{prop} \label{fHd} Assume that $0<t\leq n/2$, $H\in \cH_t$ and
$0\leq d\leq n$.  \\
(a) The degree of $f_{H,d}$ is $t$,
$\mbox{\em lm}(f_{H,d})=x_H$, and the leading coefficient is 1. \\
(b) If $D\subseteq [n]$, $|D|=d$, then $f_{H,d}(v_D)=0$.
\end{prop}

In Theorem 1.2 of \cite{HR} we determined the reduced Gr\"obner basis of the ideal $I(V{[n] \choose d})$.

\begin{thm} \label{uniform} Let $0\leq d\leq n/2$ be integers and $V:=V{[n]\choose d}$. Let $\F$ denote an arbitrary field and let $\prec$ be an arbitrary term order on the monomials of $\F[x_1,\ldots,x_n]$ for which $x_n\prec \ldots \prec x_1$. The following set of polynomials
$$
{\cal G}= \{x_2^2-x_2,\ldots ,x_n^2-x_n\}\cup \{x_J:~J\in \cD_{d}\}\cup
$$
$$
\cup\{f_{H,d}:~H\in \cH _t \mbox{ for some } 0<t\leq d\}
$$
constitutes the reduced Gr\"obner basis of the ideal $I(V)$ with respect to $\prec$. \qed
\end{thm}

Let $p$ denote an arbitrary prime.

\begin{prop} \label{bound}
Let $V:=V{[4p]\choose 2p}\subseteq \{0,1\}^{4p}\subseteq \Fp^{4p}$ and let $C\in {[4p] \choose 3p}$ be an arbitrary subset. Define $Q:=V\cup \{\underline{v}_C\}$. Let
$y\in \mbox{Sm}(Q,\prec_{deg})\setminus \mbox{Sm}(V,\prec_{deg})$. Then $deg(y)\geq p$. \qed
\end{prop}

\proof

For $0<t<p$ and $H\in \cH_t$ we define $g_H\in
\Fp[x_1,\ldots ,x_{4p}]$ as the modulo $p$ reduction of the polynomial (with
integer coefficients) $f_{H,2p}$.  By Proposition \ref{fHd} (a) the degree
of $g_H$ is $t$ and the leading term of $g_H$ is $x_H$.

By Proposition \ref{basic} and Theorem \ref{uniform} it is enough to prove that
\begin{equation} \label{egyen}
g_{H}(\underline{v}_C)=0
\end{equation}
for each $H\in \cH_t$, where $0<t<p$.

Consider the complete $p$-uniform family
\begin{equation}
\cF(p)=\{K \subseteq [4p]:~|K| \equiv 0 \mbox{ (mod }p)\}.
\end{equation}
The following Lemma follows from the Vandermonde identity
(\cite{GKP}, pp. 169-170).
\begin{lemma} \label{idez}
Let $p$ a prime. Let $x,\ j$ be integers, $0\leq j <p$.
Then
$$
{x+p \choose j} \equiv {x \choose j} \pmod p.
$$
\end{lemma} \qed

Now let $D\in \cF(p)$ and write $\underline{v}=\underline{v}_D$. Then $|D|=k'$ for some $k'$ such
that $0\leq k'\leq 4p$ and
$k'\equiv 0\mbox{ (mod}\ p)$.
We observe that $f_{H,2p}\equiv f_{H,k'}\mbox{ (mod}\ p)$, i.e., the
coefficients of the two polynomials are the same modulo $p$. Namely, for $0 \leq i \leq t$ we have
$$
{2p-i \choose t-i}\equiv {k'-i \choose t-i} \pmod p,
$$
where we used Lemma \ref{idez} and $0\leq t-i\leq p-1$.

We conclude that
$$
g_H(\underline{v})\equiv f_{H,2p}(\underline{v})\equiv f_{H,k'}(\underline{v})=0 \pmod p.
$$
Here the last equality is a consequence of Lemma \ref{fHd} (b).
Since $C\in \cF(p)$, therefore $g_H(\underline{v}_C)=0$, which was to be proved. \qed

\medskip

\section{Proofs}

{\bf Proof of Theorem \ref{main}:} Let $A_1,\ldots ,A_{m(p)}\subseteq {[4p]\choose 2p}$ denote the subsets of $[4p]$ such that for any subset $B\in {[4p] \choose 2p}$ there exists at least one $i$, $1\leq i\leq m(p)$ with $|A_i\cap B|=p$. We denote by $\underline{v_C}$ the characteristic vector of an arbitrary set $C\subseteq [4p]$. Let $\underline{v}_i:=\underline{v}_{A_i}$. Consider the following polynomial:
$$
F(x_1,\ldots,x_{4p}):=\prod_{i=1}^{m(p)} \underline{x}\cdot \underline{v}_i\in \Fp[x_1,\ldots,x_{4p}],
$$
where the central dot denotes the usual scalar product of
$\underline{x}=(x_1,\ldots ,x_n)$ and $\underline{v}_i$.

If $B\in {[4p] \choose 2p}$ is an arbitrary subset, then the previous property of the sets $A_1,\ldots, A_{m(p)}$ implies that
\begin{equation}
F(\underline{v_B})=\prod_{i=1}^{m(p)} {\underline{v}}_B\cdot {\underline{v}}_i=\prod_{i=1}^{m(p)} |A_i\cap B|\equiv \prod_{i=1}^{m(p)} |A_i\cap B|-p=0 \pmod p.
\end{equation}

\begin{prop} \label{letez}
There exists a subset $C\in {[4p] \choose 3p}$ such that
\begin{equation} \label{felt}
|C\cap A_i|\not\equiv 0 \pmod p
\end{equation}
for each $1\leq i\leq m(p)$.
\end{prop}

\proof

Let $1\leq i\leq m(p)$ be a fixed index and consider the set system
$$
\cT_i:=\{ T\in {[4p] \choose 3p}:~ |T\cap A_i|\equiv 0 \pmod p\}.
$$
Clearly it is enough to prove that
\begin{equation} \label{unio}
|\cup_{i=1}^{m(p)} \cT_i|<{4p \choose p},
\end{equation}
because then any subset from ${[4p] \choose 3p}\setminus \cup_{i=1}^{m(p)} \cT_i$
satisfies the condition (\ref{felt}). But
$$
|\cup_{i=1}^{m(p)} \cT_i|\leq \sum_{i=1}^{m(p)} |\cT_i|\leq m(p)\cdot \max_i |\cT_i|\leq 2p \max_i |\cT_i|,
$$
because $m(p)\leq 2p$.

It is easy to verify that
\begin{equation}
\{T\in {[4p]\choose 3p}:~ |T\cap A_i|=p\}\cup \{T\in {[4p]\choose 3p}:~ |T\cap A_i|=2p\}
\end{equation}
gives a disjoint decomposition of the set $\cT_i$. Since $A_i\in {[4p] \choose 2p}$ for each $1\leq i\leq m(p)$, hence
\begin{equation}
|\{T\in {[4p]\choose 3p}:~ |T\cap A_i|=p\}|=|\{T\in {[4p]\choose 3p}:~ |T\cap A_i|=2p\}|={2p\choose p}.
\end{equation}

Therefore $|\cT_i|=2\cdot {2p \choose p}$ for each $1\leq i\leq m(p)$. This implies that
\begin{equation}
2p \max_i |\cT_i|=4p{2p \choose p}<{4p \choose p},
\end{equation}
if $p>3$. \qed

Using Proposition \ref{letez}, let $C\in {[4p]\choose 3p}$ denote a fixed subset such that $|C\cap A_i|\not\equiv 0 \pmod p$ for each $1\leq i\leq m(p)$. Then clearly
\begin{equation}
F(\underline{v_C})=\prod_{i=1}^{m(p)} \underline{v}_C\cdot \underline{v}_i=\prod_{i=1}^{m(p)} |A_i\cap C|\not\equiv 0 \pmod p.
\end{equation}

Apply Theorem \ref{method} with the choices $\cF:=V{[4p] \choose 2p}\subseteq \Fp^{4p}$ and $\underline{h}:=\underline{v}_C\in \Fp^{4p}$.

Define $\cT:=V{[4p] \choose 2p}\cup {\underline{v}_C}$ and let
$$
y\in \mbox{Sm}(\cT,\prec_{deg})\setminus \mbox{Sm}(V{[4p] \choose 2p},\prec_{deg})
$$
denote the unique monomial from this difference. We proved in Theorem \ref{method} that $\mbox{deg}(F)\geq \mbox{deg}(y)$. Then $\mbox{deg}(y)\geq p$ follows from Proposition \ref{bound}. This means that $m(p)\geq \mbox{deg}(F)\geq p$, which was to be proved. \qed

{\bf Proof of Theorem \ref{main2}:} Let $\omega_0:=1$. Denote by
\begin{equation}
B_i:=\{\underline{x}=(x_1,\ldots,x_n)\in B:~ x_1\cdot \ldots \cdot x_n=\omega_i\}\subseteq B
\end{equation}
for each $0\leq i\leq p-1$. Let $B_p:=B_0$.

Let $T\subseteq B$ stand for an arbitrary set of vectors of $B$ such that for every vector $\underline{u}\in B$ there exists a $\underline{t}\in T$, with $\underline{u}\cdot \underline{t}=0$.

We must show that $|T|\geq n(p-1)$.
Define $T_i:=T\cap B_i$ for $0\leq i\leq p-1$, then clearly
$$
T=T_0\cup\ldots \cup T_{p-1}
$$
gives a disjoint decomposition of the set $T$.

Consider the following polynomials in $\underline{x}=(x_1,\ldots,x_n)$:
$$
P_i(x_1,\ldots,x_n):=\prod_{\underline{v}=(v_1,\ldots,v_n)\in T_i} (\sum_{i=1}^n v_i x_i) \in \C[x_1,\ldots ,x_n],
$$
for each $0\leq i\leq p-1$.

Then clearly $\mbox{deg}(P_i)\leq |T_i|$, therefore it is enough to prove that $\mbox{deg}(P_i)\geq \frac{n(p-1)}{p}$, because then $|T|=\sum_{i=0}^{p-1} |T_i|\geq n(p-1)$.

\begin{lemma} \label{ortog}
Let $\underline{y},\underline{z}\in B$ be arbitrary vectors. If $\underline{y}\cdot \underline{z}=0$ and $\underline{y}\in B_i$, then $\underline{z}\in B_{p-i}$ for each $0\leq i\leq p-1$.
\end{lemma}

\proof Let $\underline{y}=(y_1,\ldots,y_n)$ and $\underline{z}=(z_1,\ldots,z_n)$. Then the numbers $y_1z_1,\ldots ,y_n z_n$ are $p^{th}$ roots of unity. Suppose that these numbers give a corresponding permutation of $\lambda_0$ $\omega_0$'s
, ..., $\lambda_{p-1}$ $\omega_{p-1}$'s. Then
$$
\sum_{i=1}^n y_i z_i=\lambda_0 \omega_0+\ldots + \lambda_{p-1} \omega_{p-1}=0$$
and since $\sum_{i=0}^{p-1} \omega_i=0$, we get
$$
(\lambda_0-\lambda_{p-1})\omega_0 +\ldots +(\lambda_{p-2}-\lambda_{p-1})\omega_{p-2}=0.
$$
Indirectly, suppose that there exists an $0\leq i\leq p-2$ such that $\lambda_i\neq \lambda_{p-1}$. This means that there exists a polynomial $f\in \Q[y]$ such that $f(\omega_1)=0$ and $\mbox{deg}(f)\leq p-2$, which gives a contradiction.

Therefore $\lambda_0=\ldots =\lambda_{p-1}=\frac{n}{p}$. Consider the product $A:= \prod_{i=1}^n (y_i\cdot z_i)$. The previous argument gives that $A=(1\cdot \omega_1\cdot \ldots \cdot\omega_{p-1})^{\frac{n}{p}}=1$ and $A=\prod_{i=1}^n y_i\cdot \prod_{i=1}^n z_i=\omega_i\cdot\prod_{i=1}^n z_i$, because $\underline{y}\in B_i$. Consequently $\underline{z}\in  B_{p-i}$ \qed

We prove that $P_i(\underline{z})=0$ for every $\underline{z}\in B_{p-i}$.

Let $\underline{z}\in B_{p-i}\subseteq B$ be an arbitrary vector. Then there exists a $\underline{t}\in T\subseteq B$ such that $\underline{z}\cdot \underline{t}=0$. But Lemma \ref{ortog} implies that $\underline{t}\in B_i$. Hence $\underline{t}\in B_i\cap T=T_i$, which means that $P_i(\underline{z})=\prod_{\underline{v}\in T_i} (\underline{v}\cdot \underline{z})=0$.

Let $0\leq j\leq p-1$, $j\neq p-i$ be an arbitrary, but fixed index and let $\underline{q}\in B_j$ be an arbitrary vector. Then $P_i(\underline{q})\neq 0$, because $\underline{t}\cdot \underline{q}\neq 0$ for every $\underline{t}\in T_i=B_i\cap T$ by Lemma \ref{ortog}.

Apply Theorem \ref{method} with the choices $\cF:=B_{p-i}\subseteq \C^n$ and $\underline{h}:=\underline{q}$. Define $\cT:=B_{p-i}\cup \{\underline{q}\}\subseteq \C^n$. Consider the monomial
$$
y\in \mbox{Sm}(\cT,\prec_{deg})\setminus \mbox{Sm}(B_{p-i},\prec_{deg}).
$$

We proved in Theorem \ref{method} that $\mbox{deg}(P_i)\geq \mbox{deg}(y)$. By Corollary \ref{becs} $\mbox{deg}(y)\geq \frac{n(p-1)}{p}$, i.e., $\mbox{deg}(P_i)\geq \frac{n(p-1)}{p}$, which was to be proved. \qed

\medskip
\noindent
{\bf Acknowledgements.} The author is grateful to Lajos R\'onyai for his useful remarks.

\end{document}